\documentclass{amsart}
\usepackage{shortsalch, xy}
\xyoption{all}
\title{A recognition principle for the existence of descent data.}
\date{July 2014.}

\begin{document}

\begin{abstract}
Suppose $R\rightarrow S$ is a faithfully flat ring map. Given an $S$-module $N$, does there exists some $R$-module $M$ such that $S\otimes_R M\cong N$?
In this paper we work out (as a special case of a more general question about extensions of comonads) a criterion for the existence of such an $R$-module $M$, under
some reasonable hypotheses on the map $R\rightarrow S$. 
\end{abstract}

\maketitle
\tableofcontents

\section{Introduction.}
 Suppose $R\rightarrow S$ is an effective descent morphism of rings, for example, a faithfully
flat morphism. 
We will write $\Rep(R),\Rep(S)$ for the representation semirings of $R$ and $S$, that is, 
the isomorphism classes of finitely generated $R$-modules and $S$-modules, respectively. Here addition is given by direct sum and multiplication by
tensor product. We have a base-change (``tensoring-up'') map 
\[ \Rep(R) \stackrel{f}{\longrightarrow} \Rep(S).\]
This map may fail to be injective, but we have excellent control over its failure
to be injective. If $N\in \im f$, then descent theory identifies the preimage $f^{-1}(N)$ with the set of $S/R$-descent data that $N$ admits.
When this set is nonempty, there are useful and classical cohomological descriptions of it. 
See \cite{MR547117} for a nice exposition
of some results of this kind. A modern, very general version is in Mesablishvili's paper \cite{mesablishvili}.

If one wants to understand the map $\Rep(R)\stackrel{f}{\longrightarrow}\Rep(S)$, however, something is missing from this picture:
one needs to get some control over the failure of $f$ to be {\em surjective.}
In other words, 
we do not know how to recognize which elements of $\Rep(S)$ are indeed in the image of $f$. 
Another way of putting it is that
we want to know, given a finitely generated $S$-module $N$, whether there exists {\em any $S/R$-descent datum on $N$ at all.}
Equivalently, we want to have a simple criterion for determining whether $N\cong S\otimes_R M$ for some $R$-module $M$.
Such a recognition principle, along with the descent theory described
above, is what one needs in order to understand the relationship between $\Rep(R)$ and
$\Rep(S)$, or more generally, to understand how the module theory of a ring changes under faithfully flat extension of that ring.

The purpose of this note is to describe and prove such a recognition principle
(Theorem~\ref{main technical thm}). Our recognition principle is an abstract
statement about extensions of comonads, and as such, it has
sufficient generality to be applied to many nonclassical situations (e.g. cases in which the rings $R,S$ are not commutative, and 
may have gradings that we insist the modules respect; and there may perhaps be interesting applications which are not of an algebraic nature at all). 
A short list of the easiest cases to see that this recognition principle 
applies in is Theorem~\ref{special cases}, with consequences listed in
Corollaries~\ref{applications corollary 2}, ~\ref{applications corollary 3}, ~\ref{graded applications corollary 2}, and \ref{graded applications corollary 3}.
We also provide some explicit sample computations in Example~\ref{nonsplit example} and Example~\ref{split example}.

The most familiar setting in which our main result applies is the case in which
we have a pair of maps
\begin{equation}\label{algebra extension} A \rightarrow B\rightarrow k\otimes_A B
\end{equation}
of finite-dimensional commutative augmented algebras over a field $k$,
such that the map $A\rightarrow B$ is a split monomorphism as a map
of $A$-modules, and such that, whenever $a\in A$ is in the kernel of the augmentation map $A \rightarrow k$, then $ba\in A$ for all $b\in B$.
Then our Corollary~\ref{applications corollary 2} states that
a $B$-module $M$ is of the form $B\otimes_A N$ for some $A$-module $N$ if 
and only if $(k\otimes_A B)\otimes_B M \cong k\otimes_B M$
is a free $k\otimes_A B$-module. 
In other words, we have our criterion for the existence
of a $B/A$-descent datum on a finitely-generated $B$-module $M$: such a descent datum exists if and only if $(k\otimes_A B)\otimes_B M$
is a free $k\otimes_A B$-module. 

Consequently, if every finitely generated projective $A$-module is free and
every finitely-generated projective $B$-module is free (both of which are
often satisfied in cases of interest---for example, when $A,B$ are both
Artin $k$-algebras, since they are then nilpotent extensions of $k$ and hence
local rings), then we get an exact sequence of commutative monoids
\[ \StableRep(A) \rightarrow \StableRep(B) \rightarrow \StableRep(B\otimes_A k)\rightarrow 0,\]
where we write $\StableRep$ for the commutative monoid of stable equivalence
classes of finitely generated modules. This is our Corollary~\ref{applications corollary 3}.

We use these results in our work on stable representation theory and stable algebraic $G$-theory, \cite{g-theory1}.

We note that the initial version of this paper was entirely in terms of ``extensions of abelian categories,'' and an anonymous referee suggested rephrasing and generalizing the results to be in basically the level of generality of Definition~\ref{def of pointed extension of a comonad}. 
Namely: one has two comonads $\specialb,\specialc$ on categories $\mathcal{B},\mathcal{C}$ and a commutative diagram
\[\xymatrix{ \specialc\hyphencoalg \ar[r]^F & \specialb\hyphencoalg \ar[d]^{U_{\specialb}}\ar[r]^{t^*} & \specialc\hyphencoalg\ar[d]^{U_{\specialc}} \\
 & \mathcal{B} \ar[r]^{s^*} & \mathcal{C}, }
\]
and given an object $X$ of $\mathcal{B}$ such that $s^*X$ admits the structure of
a $\specialc$-coalgebra, one wants to know whether $X$ admits the structure of a 
$\specialb$-coalgebra. Our Theorem~\ref{main technical thm} gives necessary and sufficient
conditions for this to be so, and it in turns implies our main result with
algebraic applications, Theorem~\ref{special cases}. It would be interesting to know if there exist non-algebraic applications for our Theorem~\ref{main technical thm}. 

We are grateful to the anonymous referee for their apt and helpful suggestions, and also to C. Weibel for his editorial help.

\section{The main definitions and the main result.}


\begin{definition}\label{def of pointed extension of a comonad}
Let $\mathcal{C}$ be a category, $\specialc$ a comonad on $\mathcal{C}$.
By a {\em pointed extension of $\specialc$} we mean the following data:
\begin{itemize}
\item a category $\mathcal{B}$ and a comonad $\specialb$ on $\mathcal{B}$, and
\item functors $F,t^*,s^*$ which admit right adjoints $G,t_*,s_*$, respectively,
which fit into the following diagram and make it commute:
\[\xymatrix{ \specialc\hyphencoalg \ar[r]^F & \specialb\hyphencoalg \ar[d]^{U_{\specialb}}\ar[r]^{t^*} & \specialc\hyphencoalg\ar[d]^{U_{\specialc}} \\
 & \mathcal{B} \ar[r]^{s^*} & \mathcal{C}, }
\]
where by $\specialc\hyphencoalg$ and $\specialb\hyphencoalg$ we mean the categories of $\specialc$-coalgebras and $\specialb$-coalgebras, respectively, and by $U_{\specialc},U_{\specialb}$ we mean the usual forgetful functors. (In what follows, 
we will write $W_{\specialc}, W_{\specialb}$ 
for the usual right adjoints to these functors.)
\end{itemize}
The above data is required to satisfy the following conditions:
\begin{itemize}
\item {\em (Pointedness.)} $t^* \circ F \simeq \id_{\specialc\hyphencoalg}$.
\item {\em (Beck-Chevalley condition.)} $U_{\specialb}\circ t_* = s_*\circ U_{\specialc}$.
\item {\em (Affineness.)} The functors $G$ and $W_{\specialb}$ preserve epimorphisms.
\item {\em (Semisimplicity.)} Every epimorphism in the category
of $\specialc$-coalgebras is split.
\item {\em (Nakayama axiom.)} The unit map $X \rightarrow s_*s^*X$ of the adjunction is an epimorphism for all objects $X$ of $\mathcal{B}$, and $s_*s^*$ has the property that if $s_*s^*f$ is an isomorphism, then $f$ is an epimorphism.
\end{itemize}
\end{definition}
The names of the axioms which suggest a geometric origin are because these axioms are motivated by, and satisfied in, algebraic situations of geometric interest. Specifically, in the case where we have a field $k$, an augmented commutative $k$-algebra $R$, and a monomorphism of commutative $k$-algebras
$R \rightarrow S$, we let $\mathcal{B}$ be the category of $S$-modules and we let
$\mathcal{C}$ be the category of $k\otimes_R S$-modules, with $t^*$ the
base change functor $M \mapsto (k\otimes_R S)\otimes_S M$. We let 
$\specialb$ be the base-change comonad $M \mapsto S\otimes_R M$ on 
$S$-modules, and we let $\specialc$ be the base-change comonad
$M \mapsto (k\otimes_R S)\otimes_k M$ on $k\otimes_R S$-modules.
Then the category of $\specialc$-coalgebras is equivalent to the category of $k$-modules (thanks to the anonymous referee, who pointed out this useful fact!), and
the category of $\specialb$-coalgebras is equivalent to the category of $R$-modules
if the morphism $R\rightarrow S$ is an effective descent morphism.

In that situation, the pointedness axiom in Definition~\ref{def of pointed extension of a comonad}
is satisfied because of $R$ being an augmented, i.e., pointed, $k$-algebra.
The affineness axiom is satisfied because the scheme morphisms
$\Spec S \rightarrow \Spec R$ and $\Spec R \rightarrow \Spec k$ are
affine morphisms, so their induced direct image functors on the quasicoherent
module categories are right exact, i.e., preserve epimorphisms.
The semisimplicity axiom is satisfied because the category of $k$-vector spaces
is semisimple.
On restricting to the finitely-generated module categories, 
the Nakayama condition is satisfied when the kernel of $S\rightarrow k\otimes_R S$
is contained in the Jacobson radical of $S$, by Nakayama's lemma.

Now we define what it means for an object to be ``formally supported''
by a comonad. The definition looks technical at a glance, but it has two
advantages: first, the objects $X$ formally supported by a comonad $\specialb$ have the property that
they admit the structure of a $\specialb$-coalgebra if and only if 
$s^*X$ admits the structure of a $\specialc$-coalgebra. See Theorem~\ref{main technical thm}
for this result. Second, the condition required for an object to be 
formally supported by a comonad is actually concrete enough that it is checkable
in situations of concrete interest, and frequently holds for {\em all} objects on which the comonad is defined. See our Theorem~\ref{special cases}, where we 
verify that all objects are formally supported by certain 
base-change comonads arising in a wide class
of very concrete and explicit algebraic cases of interest.

\begin{definition}\label{def of formal support}
Let $\specialc$ be a comonad and $(\specialb, F,t^*,s^*)$ a pointed extension 
of $\specialc$. 
(Our notation is as in Definition~\ref{def of pointed extension of a comonad}.)
We say that an object $M$ of $\mathcal{B}$ is {\em formally supported 
by $\specialb,F,t^*,s^*,\specialc$}, or {\em formally supported 
by $\specialb$} for short, if the following condition is satisfied:
\begin{itemize}
\item {\em (Lifting axiom.)}
For every $\specialc$-coalgebra $Y$ and every epimorphism
$\tilde{\sigma}: U_{\specialb} FY \rightarrow M$ such that $s^*\tilde{\sigma}$ is an 
isomorphism, there exists a map
$\tilde{\phi}: M \rightarrow U_{\specialb}W_{\specialb}M$
making the diagram
\[\xymatrix{ U_{\specialb} FY \ar[d]^{\tilde{\sigma}} \ar[r]^{U_{\specialb}\eta_{FY}} & U_{\specialb}W_{\specialb}U_{\specialb}FY\ar[d]^{U_{\specialb}W_{\specialb}\tilde{\sigma}} \\ 
M\ar@{-->}[r]^{\tilde{\phi}} & U_{\specialb}W_{\specialb} M }\]
commute.
\end{itemize}
\end{definition}

If enough objects are formally supported by a comonad, we say that the pointed
extension of comonads is {\em short}:
\begin{definition}
Let $\specialc$ be a comonad and $(\specialb, F,t^*,s^*)$ a pointed extension 
of $\specialc$. 
We say that {\em $\specialb$ is an short extension of $\specialc$}
if, whenever $M$ is an object of $\mathcal{B}$ 
such that $s^*M$ is isomorphic to $s^*U_{\specialb}FY$ for some $\specialc$-coalgebra $Y$,
we have that $M$ is formally supported by $\specialb$.
\end{definition}

Now here is our main technical result, which implies the concrete results in the
next section:
\begin{theorem}\label{main technical thm}
Let $\specialc$ be a comonad and $(\specialb, F,t^*,s^*)$ a pointed extension 
of $\specialc$. 
(Our notation is as in Definition~\ref{def of pointed extension of a comonad}.)
Suppose $M$ is an object of $\mathcal{B}$ which is formally supported by $\specialb$. Then $M$ admits the structure of a $\specialb$-coalgebra if and only if $s^*M$ admits the structure of a $\specialc$-coalgebra. 
\end{theorem}
\begin{proof}
Suppose $s^*M$ admits the structure of a $\specialc$-coalgebra. 
Then we can choose a $\specialc$-coalgebra $Y$ such that
$U_{\specialc} Y = s^*M$ in $\mathcal{C}$.
Let $\sigma: U_{\specialb}FY\rightarrow s_*s^*M$ be the composite
\begin{equation}\label{commutative diagram 2} \xymatrix{
U_{\specialb}FY \ar[r]^{\eta_{U_{\specialb}FY}} \ar@{-->}[d]^{\tilde{\sigma}}\ar[rd]^{\sigma} &
 s_*s^*U_{\specialb} FY \ar[r]^{\cong} & s_*U_{\specialc}t^*FY \ar[d]^{\cong} \\
M \ar[r]^{\eta_M} & s_*s^*M & s_*U_{\specialc}Y\ar[l]^{=} ,
}\end{equation}
i.e., $s = \eta_{U_{\specialb}FY}$.

We have the commutative diagram of hom-sets
\begin{equation}\label{commutative diagram 1}\xymatrix{
\hom_{\mathcal{B}}(U_{\specialb} FY, M) \ar[r]\ar[d]^{\cong} &
 \hom_{\mathcal{B}}(U_{\specialb}FY, s_*s^*M) \ar[d]^{\cong} \\
\hom_{\specialc\hyphencoalg}(Y, GW_{\specialb} M) \ar[r] &  
 \hom_{\specialc\hyphencoalg}(Y, GW_{\specialb}s_*s^*M) }\end{equation}
in which the bottom horizontal map is split epic,
since $M \rightarrow s_*s^*M$ is epic by 
the Nakayama axiom in Definition~\ref{def of pointed extension of a comonad},
since $G,W_{\specialb}$ both preserve epimorphisms by the
affineness axiom in Definition~\ref{def of pointed extension of a comonad},
and since every epimorphism in $\specialc$-coalgebras is
split by the semisimplicity axiom in Definition~\ref{def of pointed extension of a comonad}.
So the top horizontal map in diagram~\ref{commutative diagram 1} is also
a split epimorphism.
So the map $\sigma\in \hom_{\mathcal{B}}(U_{\specialb}FY, s_*s^*M)$
lifts to a map $\tilde{\sigma}\in \hom_{\mathcal{B}}(U_{\specialb}FY, M)$
which fills in the dotted arrow in diagram~\ref{commutative diagram 2}.

We claim that $\tilde{\sigma}$ is epic. By the commutativity of diagram~\ref{commutative diagram 2}, 
$s_*s^*\tilde{\sigma}$ makes the diagram
\[\xymatrix{
s_*s^*U_{\specialb}FY\ar[d]^{s_*s^*\tilde{\sigma}}\ar[r]^{\cong} & s_*U_{\specialc}t^*FY \ar[d]^{\cong} \\
s_*s^*M & s_*U_{\specialc}Y\ar[l]^{}_{=} 
}\]
commute, hence $s_*s^*\tilde{\sigma}$ is an isomorphism. Now since $s_*s^*\tilde{\sigma}$ is an isomorphism,
$\tilde{\sigma}$ is an epimorphism by the Nakayama axiom.

We now have the diagram in $\mathcal{B}$:
\[\xymatrix{
U_{\specialb}FY \ar[r]^{U_{\specialb}\eta^{\specialb}_{FY}} \ar[d]^{\tilde{\sigma}} &
 U_{\specialb}W_{\specialb}U_{\specialb}FY \ar[d]^{U_{\specialb}W_{\specialb}\tilde{\sigma}} \\
M \ar@{-->}[r]^{\tilde{\psi}} & U_{\specialb}W_{\specialb}M, }\]
and the lifting axiom in Definition~\ref{def of formal support} now guarantees
the existence of a map $\tilde{\psi}$ filling in the dotted arrow
and making the diagram commute.

We claim that $\tilde{\psi}$ is in fact a $\specialb$-coalgebra structure map
on $M$.
To prove this claim, we must check that it is counital and coassociative.

We check counitality first: we have the commutative diagram with exact columns
\[\xymatrix{
 U_{\specialb}FY \ar[r]^>>>>{\psi_{std}} \ar[d]^{\tilde{\sigma}} \ar@/^3pc/[rrr]^{\id} &
  U_{\specialb}W_{\specialb}U_{\specialb}FY\ar[d]^{U_{\specialb}W_{\specialb}\tilde{\sigma}} \ar[rr]^>>>>>>>>>>{\epsilon U_{\specialb}F Y} & &
  U_{\specialb}F \ar[d]^{\tilde{\sigma}}Y \\
 M \ar[r]^{\tilde{\psi}} & U_{\specialb}W_{\specialb}M\ar[rr]^{\epsilon M} & & M  ,
}\]
where $\epsilon$ is the counit natural transformation of the adjunction of $U_{\specialb},W_{\specialb}$.
We are also writing $\psi_{std}$ for the standard $\specialb$-coalgebra structure
map on $U_{\specialb}FY$, that is,
$\psi_{std} = U_{\specialb}\eta^{\specialb}_{FY}$.
We get the identity map across the top of the diagram because of $\psi_{std}$ being itself an $U_{\specialb}W_{\specialb}$-coalgebra structure map, hence itself counital.
From the commutativity of this diagram we get the equality
\begin{eqnarray*} \epsilon M\circ \tilde{\psi} \circ \tilde{\sigma} & = & \tilde{\sigma}\circ \epsilon U_{\specialb}FY\circ \psi_{std} \\
 & = & \tilde{\sigma}.\end{eqnarray*}
Now since $\tilde{\sigma}$ is epic, i.e., right-cancellable, we have that $\epsilon M\circ \tilde{\psi} = \id$, which is precisely the statement of counitality
for $\tilde{\psi}$.

Now we check coassociativity. Due to basic properties of adjunctions and their comonads and also coassociativity of $\psi_{std}$ since it itself
 is a $\specialb$-coalgebra structure map, we have the equalities
\begin{eqnarray*}
U_{\specialb}W_{\specialb}\tilde{\psi} \circ \tilde{\psi} \circ \tilde{\sigma} 
& = & U_{\specialb}W_{\specialb}\tilde{\psi} \circ U_{\specialb}W_{\specialb}\tilde{\sigma} \circ \psi_{std} \\
& = & U_{\specialb}W_{\specialb}(\tilde{\psi}\circ \tilde{\sigma}) \circ \psi_{std} \\
& = & U_{\specialb}W_{\specialb}(U_{\specialb}W_{\specialb}\tilde{\sigma} \circ \psi_{std})\circ\psi_{std} \\
& = & U_{\specialb}W_{\specialb}U_{\specialb}W_{\specialb}\tilde{\sigma} \circ 
         U_{\specialb}\eta W_{\specialb}U_{\specialb}F Y\circ \psi_{std} \\
& = & U_{\specialb}\eta W_{\specialb} M \circ U_{\specialb}W_{\specialb}\tilde{\sigma} \circ \psi_{std} \\
& = & U_{\specialb}\eta W_{\specialb} M \circ \tilde{\psi}\circ \tilde{\sigma},\end{eqnarray*}
and since $\tilde{\sigma}$ is epic, i.e., right-cancellable, this tells us that
\[ U_{\specialb}W_{\specialb}\tilde{\psi} \circ \tilde{\psi}
 = U_{\specialb}\eta W_{\specialb} M \circ \tilde{\psi},\]
which is precisely the statement that $\tilde{\psi}$ is coassociative.

Conversely, if $M$ is a $\specialb$-coalgebra, then $M$ is in the image of $U_{\specialb}$. Choose $N$ in the category of $\specialb$-coalgebras such that
$M = U_{\specialb}(N)$. Then we have
\begin{eqnarray*}
 s^*M & = & s^*U_{\specialb} N \\
      & = & U_{\specialc}t^* N,\end{eqnarray*}
so $s^*M$ is in the image of $U_{\specialc}$.
So $s^*M$ admits the structure of a $\specialc$-coalgebra.
\end{proof}

\begin{corollary}\label{corollary of having an short extension}
Let $\specialc$ be a comonad and let $\specialb$ be an short extension of
$\specialc$. Then an object $M$ of $\mathcal{B}$ admits the structure of 
a $\specialb$-coalgebra if and only if $s^*M$ admits the structure of a 
$\specialc$-coalgebra.
\end{corollary}

Corollary~\ref{corollary of having an short extension} 
makes it clear that the relationship between $\specialc$-coalgebras
and $\specialb$-coalgebras is very simple when $\specialb$ is an short 
extension of $\specialc$.
We claim that short extensions of
comonads are actually not unusual, in fact
occuring in practical, concrete situations.
See Theorem~\ref{special cases} for a wide class of examples,
and Examples~\ref{nonsplit example} and~\ref{split example} for
a couple of completely explicit computations.

\section{Special cases and applications.}

We introduce a quick definition of a certain class of monoids which are suitable for being the monoids of grading for graded objects, e.g.
$\mathbb{N}$ and $\mathbb{Z}$.

\begin{definition} We will say that a commutative monoid $\mathbb{M}$ is {\em finitely-generated and weakly free,} or {\em FGWF} for short,
if $\mathbb{M}$ is isomorphic to a finite Cartesian product of copies of $\mathbb{N}$ and $\mathbb{Z}$.\end{definition}

Now we introduce a useful definition which plays an essential role in Theorem~\ref{special cases}. 
\begin{definition}\label{def of short map}
Suppose $k$ is a field, and $A,B$ are $k$-algebras.
Suppose $A$ is equipped with a $k$-algebra 
augmentation map $\epsilon: A\rightarrow k$.
We say that a monomorphism $f: A\rightarrow B$
of $k$-algebras is {\em (right) short} if,
for each element $b\in B$ and each element
$a\in \ker \epsilon$, the element $ba\in B$
is in $\im f\subseteq B$.
\end{definition}
One can define a ``left short'' morphism by simply swapping right-multiplication
out for left-multiplication in Definition~\ref{def of short map}.

Note that any short morphism $R \rightarrow S$ has the property that
$S\rightarrow S\otimes_R k$ is a (right) square-zero extension of algebras,
but the converse is not necessarily true: if $S\rightarrow S\otimes_R k$
is a right square-zero extension, it is not always true that $R \rightarrow S$
is a right short morphism, even when $R\rightarrow S$ is injective. For example,
let $S = k[x]/x^4$ and let $R = k[y]/y^2$, and let $R\rightarrow S$ send
$y$ to $x^2$. This map is not short, but $S \rightarrow S\otimes_R k$ is
a square-zero extension. On the other hand, the map of $k$-algebras
$k[y]/y^2 \rightarrow k[x]/x^3$ sending $y$ to $x^2$ is a good example of a short morphism.

Here are some examples of pointed extensions of comonads, to which Theorem~\ref{main technical thm} applies.
\begin{theorem}\label{special cases}
Suppose $k$ is a field. 
\begin{itemize}
\item \label{extension of rings}
Suppose $A$ is an augmented commutative algebra over $k$, and 
$f: A\rightarrow B$ is an effective descent morphism 
of commutative $k$-algebras. Suppose both $A$ and $B$ are finite-dimensional
as $k$-vector spaces. Let $C$ be the algebra
$B\otimes_A k$, and suppose the kernel of the ring map
$B\rightarrow B\otimes_A k$ is contained in the Jacobson radical of $B$.
Then the commutative diagram
of categories and functors
\[\xymatrix{
\fgMod(k) \ar[r]^F & \fgMod(A)\ar[d]^{U_{\specialb}} \ar[r]^{t^*} & \fgMod(k) \ar[d]^{U_{\specialc}} \\
& \fgMod(B)\ar[r]^{s^*} & \fgMod(C) }\]
expresses $\specialb$ as a pointed extension
of $\specialc$,
where we are writing $\specialb$ for the base change
comonad $M \mapsto \specialb M = M\otimes_A B$
on $\fgMod(B)$, 
and $\specialc$ for the base-change comonad
$M \mapsto \specialc M = M\otimes_B C$
on $\fgMod(C)$.
Here $F,U_{\specialb},U_{\specialc},t^*$ are all the
obvious base-change (tensoring-up) functors.

Furthermore, if $f$ is short, then
$\specialb$ is an short extension of
$\specialc$.
\item \label{extension of graded rings}
The previous example works as well in the graded setting.
Suppose $\mathbb{M}$ is an FGWF monoid (e.g. $\mathbb{M} = \mathbb{N}$ or $\mathbb{M} = \mathbb{Z}$)
and
suppose that
$A$ is an $\mathbb{M}$-graded augmented algebra over $k$ and suppose 
$f: A\rightarrow B$ is an $\mathbb{M}$-grading-preserving effective descent morphism (for right modules) of 
$\mathbb{M}$-graded $k$-algebras. Suppose both $A$ and $B$ are finite-dimensional
as $k$-vector spaces. Let $C$ be the $\mathbb{M}$-graded algebra
$B\otimes_A k$, and suppose the kernel of the ring map
$B\rightarrow B\otimes_A k$ is contained in the Jacobson radical of $B$.
Then
\[\xymatrix{
\mathbb{M}-gr\fgMod(k) \ar[r]^F & \mathbb{M}-gr\fgMod(A)\ar[d]^{U_{\specialb}} \ar[r]^{t^*} & \mathbb{M}-gr\fgMod(k) \ar[d]^{U_{\specialc}} \\
& \mathbb{M}-gr\fgMod(B)\ar[r]^{s^*} & \mathbb{M}-gr\fgMod(C) }\]
expresses $\specialb$ as a pointed extension
of $\specialc$,
where we are writing $\specialb$ for the base change
comonad $M \mapsto \specialb M = M\otimes_A B$
on $\mathbb{M}-gr\fgMod(B)$, 
and $\specialc$ for the base-change comonad
$M \mapsto \specialc M = M\otimes_B C$
on $\mathbb{M}-gr\fgMod(C)$.
Here $F,U_{\specialb},U_{\specialc},t^*$ are all the
obvious base-change (tensoring-up) functors.

Furthermore, if $f$ is short, then
$\specialb$ is an short extension of
$\specialc$.
\end{itemize}
\end{theorem}
\begin{proof}
Suppose $A$ is an augmented commutative algebra over $k$, and 
$f: A\rightarrow B$ is an effective descent morphism of commutative $k$-algebras. Suppose both $A$ and $B$ are finite-dimensional
as $k$-vector spaces. Let $C$ be the algebra
$B\otimes_A k$.
Then the category of $\specialc$-coalgebras is naturally equivalent to the category of $k$-vector spaces, since
the unit $k$-algebra map $k\rightarrow C$
is clearly faithfully flat, hence an effective descent morphism.
Furthermore, the category of $\specialb$-coalgebras is naturally equivalent to the category of $A$-modules, by the assumption that $f$ is an effective descent morphism. It is trivial that 
these statements remain true with
the adjective ``finitely generated'' included throughout.

Now we check the relevant conditions, from Definition~\ref{def of pointed extension of a comonad},
for $\specialb$ to be a pointed extension of $\specialc$.
\begin{itemize}
\item Pointedness is immediate from the composite of
the $k$-algebra maps $k \rightarrow C \rightarrow k$
being the identity on $k$.
\item The Beck-Chevalley condition is a classical
property of restriction and induction of scalars.
\item The affineness condition is 
classical: the restriction of scalars
functor induced by a map of rings (or $k$-algebras)
is exact.
\item Semisimplicity is due to 
$k$ being a field, hence the category of $k$-modules (or of finitely-generated $k$-modules)
is semisimple.
\item For the Nakayama condition: since
the morphism of $k$-algebras $B \rightarrow C$
is surjective, we have that
$M \cong B\otimes_B M \rightarrow C\otimes_B M
\cong s_*s^*M$ is surjective for all $B$-modules 
$M$. So the unit map $\id\rightarrow s_*s^*$ is always an epimorphism, satisfying the
first part of the Nakayama condition.

For the second part, we must use the assumption that all our modules are finitely generated.
Now the functor $s_*s^*$ is the functor
\[  M\mapsto M\otimes_B (B\otimes_{A} k)\cong M\otimes_A k\cong M/M(\ker f).\]
If $M/M(\ker g)\cong 0$, then the inclusion $M(\ker g)\hookrightarrow M$ is an isomorphism. Nakayama's Lemma
now applies: since $(\ker g)$ is contained in the Jacobson radical of $B$, $(\ker g)M = M$ and $M$ finitely generated together imply that $M = 0$.
\end{itemize}

Now suppose that $f$ is short,
and suppose that $M$ is a finitely-generated
$B$-module equipped with an epimorphism
$U_{\specialb}FY = Y\otimes_k B \stackrel{i}{\longrightarrow} M$
of $B$-modules such that $s^*i = i\otimes_B C$
is an isomorphism, where $Y$ is a 
finitely-generated $k$-module.
We will write $\epsilon: A \rightarrow k$
for the augmentation map on $A$.
Then $\ker i$ consists entirely of elements of
the free $B$-module $Y\otimes_k B$
which are divisible by elements in $\ker \epsilon$,
that is, $\ker i$ consists entirely of 
elements of the form 
$ya$ for some $y\in Y\otimes_k B$ and 
some $a\in\ker\epsilon$.

Now we choose a basis $y_1,\dots ,y_n$ for $Y$.
If $ya$ is some element of $\ker i$, 
we write it in the form
\[ ya = \left( \sum_{j=1}^n \alpha_j y_j\right)a\]
and we observe that, since $f$ is short,
each element $\alpha_j y_ja\in B\{ y_j\}\cong B$
is contained in $\im f$.
Hence $ya$ is contained in 
\[ Y\otimes_k \im f : Y\otimes_k A\rightarrow
 Y\otimes_k B.\]

Hence, in the commutative diagram
\begin{equation}\label{comm diag 10} \xymatrix{
 & \ker i \ar[d]\ar@{-->}[dl] & \\
Y\otimes_k A \ar[r] & Y\otimes_k B\ar[r] \ar[d]^i &
 Y\otimes_k B\otimes_A B\ar[d] \\
 & M & M\otimes_A B }\end{equation}
there exists a map as in the dotted arrow
which makes the diagram commute. Since the 
row in diagram~\ref{comm diag 10}
is exact and since the middle column
in diagram~\ref{comm diag 10}
is a short exact sequence,
the composite map $\ker i \rightarrow M\otimes_A B$
is zero,
hence there exists a map
$M \rightarrow M\otimes_A B = s_*s^*M$
fitting into the bottom row of diagram~\ref{comm diag 10} and making the diagram commute.
This is precisely
the lifting axiom in Definition~\ref{def of formal support}.
So if $M \otimes_B C$ is a free $B$-module,
then $M$ is formally supported by $\specialb$.
So $\specialb$ is an short extension
of $\specialc$.
\end{proof}

It is convenient to have a characterization of effective descent morphisms which is applicable
to short $k$-algebra morphisms, so that one can use Theorem~\ref{special cases} in practical situations.
One knows, from classical descent theory, that faithfully flat ring morphisms are effective descent 
morphisms. This is not so helpful, however, since it is very rare for a short $k$-algebra map to be flat.
Instead, the more general Joyal-Tierney theorem is quite helpful for determining when a short morphism of commutative $k$-algebras
is an effective descent morphism.
See \cite{MR1742958} for the following result, which Mesablishvili attributes
to Joyal and Tierney, but which was apparently never written up by Joyal and Tierney themselves:
\begin{theorem} \label{joyal-tierney thm} {\bf (Joyal-Tierney.)} Suppose $f: R \rightarrow S$ is a homomorphism of commutative rings.
Then $f$ is an effective descent morphism if and only if $f$ is a pure monomorphism.
\end{theorem}
We remind the reader that a monomorphism of $R$-modules $f: M \rightarrow N$ is said to be {\em pure} if,
for every $R$-module $M_0$, the tensor product $M_0\otimes_R f: M_0 \otimes_R M\rightarrow M_0\otimes_R N$
is a monomorphism. So the condition of Theorem~\ref{joyal-tierney thm} is that the map $R\rightarrow S$
is a pure monomorphism, as a map of $R$-modules.

The following is a corollary of the Joyal-Tierney theorem but is also not hard to prove directly.
\begin{corollary}
Suppose $f: R\rightarrow S$ is a homomorphism of commutative rings which is a split monomorphism as a map of $R$-modules.
Then $f$ is an effective descent morphism.
\end{corollary}
\begin{proof}
If $f$ is a split monomorphism as a map of $R$-modules, we can write $S$ as $R \oplus M$ for some $R$-module $M$,
with the map $f$ being the inclusion of the summand $R$. Then, for any $R$-module $M_0$, the map
\[ M_0 \cong M_0\otimes_R R \stackrel{M_0\otimes_R f}{\longrightarrow} M_0\otimes_R (R\oplus M)\cong M\oplus M_0\otimes_R M \]
is also inclusion of a summand, hence still a monomorphism. Hence $f$ is a pure monomorphism, hence an effective descent morphism
by Theorem~\ref{joyal-tierney thm}.
\end{proof}

Now we can list some more corollaries of Theorem~\ref{special cases}.

\begin{corollary}\label{applications corollary 2}
Suppose, for any ring $R$, we write $\Rep(R)$ for the commutative monoid of isomorphism classes of finitely-generated
right $R$-modules, with addition given by direct sum.
Suppose $k$ is a field, $A$ is an augmented commutative $k$-algebra, and $B$ is a commutative $k$-algebra
equipped with a short $k$-algebra monomorphism $f: A\rightarrow B$. Suppose $A,B$ are both finite-dimensional as $k$-vector spaces,
and suppose that $f$, when regarded as a map of $A$-modules, is a split monomorphism.
Finally,
suppose the kernel of the ring map
$B\rightarrow B\otimes_A k$ is contained in the Jacobson radical of $B$.
Then the image of the base-change (``tensoring up'') map of monoids $\Rep(A)\rightarrow \Rep(B)$
consists of exactly the isomorphism classes of $B$-modules $M$ such that $M\otimes_A k$ is a free 
$B\otimes_A k$-module.
\end{corollary}

\begin{corollary}\label{applications corollary 3}
Suppose $k,A,B,f$ are as in Corollary~\ref{applications corollary 2}.
Write $\StableRep(A)$ (resp. $\StableRep(B)$) for the stable representation monoid of $A$ (resp. $B$), that is,
the monoid of stable equivalence classes of finitely generated $A$-modules (resp. finitely generated $B$-modules).
Suppose every finitely generated projective $B\otimes_A k$-module is free. Then
the sequence of monoid maps
\begin{equation}\label{exact seq 1} \StableRep(A) \rightarrow \StableRep(B) \rightarrow \StableRep(B\otimes_A k)\rightarrow 0\end{equation}
is exact.
\end{corollary}

\begin{warning}
The reader should be careful to note that, while the sequence~\ref{exact seq 1} is an exact sequence of commutative monoids---that is, the image of each map
is equal to the kernel of the next---the same is not necessarily true after one takes the Grothendieck group completion of those monoids.
Furthermore, exact sequences of commutative monoids really are quite different from exact sequences of abelian groups. In particular,
$\StableRep(B)$ is capable of being non-finitely-generated even when both $\StableRep(A)$ and $\StableRep(B\otimes_A k)$ are finitely generated. See Example~\ref{split example}
for an example of this.
\end{warning}

We also have the graded analogues of Corollary~\ref{applications corollary 2} and Corollary~\ref{applications corollary 3}.

\begin{corollary}\label{graded applications corollary 2}
Suppose $\mathbb{M}$ is an FGWF monoid (e.g. $\mathbb{M} = \mathbb{N}$ or $\mathbb{M} = \mathbb{Z}$).
Suppose, for any $\mathbb{M}$-graded ring $R$, we write $\Rep(R)$ for the commutative monoid of isomorphism classes of finitely-generated $\mathbb{M}$-graded
right $R$-modules, with addition given by direct sum.
(This commutative monoid has a natural suspension action of $\mathbb{N}[\mathbb{N}]$, making it a $\mathbb{N}[\mathbb{M}]$-semimodule.)
Suppose $k$ is a field, $A$ is an $\mathbb{M}$-graded augmented commutative $k$-algebra, and $B$ is an $\mathbb{M}$-graded commutative $k$-algebra
equipped with a grading-preserving short $k$-algebra monomorphism $f: A\rightarrow B$. Suppose $A,B$ are both finite-dimensional as $k$-vector spaces,
and suppose that $f$, when regarded as a map of $A$-modules, is a split monomorphism.
Then the image of the base-change (``tensoring up'') map of monoids $\Rep(A)\rightarrow \Rep(B)$
consists of exactly the isomorphism classes of $B$-modules $M$ such that $M\otimes_A k$ is a free 
$B\otimes_A k$-module.
\end{corollary}

\begin{corollary}\label{graded applications corollary 3}
Suppose $\mathbb{M},k,A,B,f$ are as in Corollary~\ref{graded applications corollary 2}.
Write $\StableRep(A)$ (resp. $\StableRep(B)$) for the stable representation monoid of $A$ (resp. $B$), that is,
the monoid of stable equivalence classes of finitely generated $\mathbb{M}$-graded $A$-modules (resp. finitely generated $B$-modules).
Suppose every finitely generated $\mathbb{M}$-graded projective $B\otimes_A k$-module is free. Then
the sequence of monoid maps
\begin{equation*} \StableRep(A) \rightarrow \StableRep(B) \rightarrow \StableRep(B\otimes_A k)\rightarrow 0\end{equation*}
is exact.
\end{corollary}

We end with a couple of explicit computations.
\begin{example} \label{nonsplit example}
Let $k$ be a field, and let $f: k[y]/y^2\rightarrow k[x]/x^3$
be the $k$-algebra map sending $y$ to $x^2$. Then it is easily seen that:
\begin{itemize}
\item $f$, regarded as a map of $k[y]/y^2$-modules, is a split monomorphism, 
\item $k[y]/y^2$ and $k[x]/x^3$ are each finite-dimensional as $k$-vector spaces, 
\item $f$ is a short morphism, and
\item $(x^2)$ is contained in the Jacobson radical of $k[x]/x^3$.
\end{itemize}
Hence Corollary~\ref{applications corollary 2} applies.
By the classification of modules over a principal ideal domain, we have:
\begin{eqnarray*}
\Rep(k[y]/y^2) & = & \mathbb{N}\{ k,k[y]/y^2\} ,\\
\Rep(k[x]/x^3) & = & \mathbb{N}\{ k,k[x]/x^2, k[x]/x^3\}, \\
\Rep(k[x]/x^2) & = & \mathbb{N}\{ k,k[x]/x^2\}, \end{eqnarray*}
and the base-change maps are easily computed:
\[\xymatrix{ \Rep(k[y]/y^2) \ar[r] & \Rep(k[x]/x^3)\ar[r] & \Rep(k[x]/x^2)  \\
 & k\ar@{|->}[r] & k \\
k \ar@{|->}[r] & k[x]/x^2 \ar@{|->}[r] & k[x]/x^2 \\
k[y]/y^2 \ar@{|->}[r] & k[x]/x^3 \ar@{|->}[r] & k[x]/x^2. }\]
This verifies Corollary~\ref{applications corollary 2} directly, in this special case.

Every finitely generated projective module over these algebras is free, hence on passing to the stable
representation monoids, we get the description:
\[\xymatrix{ \StabRep(k[y]/y^2)\ar[d]^{\cong} \ar[r] & \StabRep(k[x]/x^3)\ar[d]^{\cong}\ar[r] & \StabRep(k[x]/x^2)\ar[d]^{\cong}  \\
\mathbb{N}\{ k\} \ar[r] & \mathbb{N}\{ k,k[x]/x^2\} \ar[r] & \mathbb{N}\{ k\} \\
 & k\ar@{|->}[r] & k \\
k \ar@{|->}[r] & k[x]/x^2 \ar@{|->}[r] & 0 . }\]
This verifies Corollary~\ref{applications corollary 3} directly, in this special case.
\end{example}

\begin{example} \label{split example}
Let $k$ be a field, and let $f: k[y]/y^2\rightarrow k[x,y]/x^2,xy,y^2$
be the $k$-algebra map sending $y$ to $y$. We give these algebras a $\mathbb{Z}$-grading
by letting $x,y$ be in any two distinct positive degrees $|x|,|y|$. 
Then it is easily seen that:
\begin{itemize}
\item $f$, regarded as a map of $k[y]/y^2$-modules, is a split monomorphism, 
\item $k[y]/y^2$ and $k[x,y]/x^2,xy,y^2$ are each finite-dimensional as $k$-vector spaces, 
\item $f$ is a short morphism, and
\item $(y)$ is contained in the Jacobson radical of $k[x,y]/x^2,xy,y^2$.
\end{itemize}
Hence Corollary~\ref{graded applications corollary 2} applies.

We use Margolis' classification of finitely-generated graded
modules over a graded exterior algebra
with generators in distinct positive degrees, from \cite{MR738973}.
Following Margolis, if $n\in\mathbb{N}$ and $\epsilon_0,\epsilon_1\in\{0,1\}$, 
we write $L(n,\epsilon_0,\epsilon_1)$ for the $k[x,y]/x^2,xy,y^2$-module
with $n+1$ generators $x_0, \dots ,x_n$ and with relations
\[ yx_i = xx_{i+1} \]
for $i=0, \dots ,n-1$, and an additional relation $xx_0$ if $\epsilon_0 = 0$, and
$yx_n$ if $\epsilon_1 = 0$. Margolis proves that, over $k[x,y]/x^2,y^2$, 
every finitely-generated graded module decomposes uniquely into a direct sum of these
(colloquially, ``lightning flash'') modules $L(n,\epsilon_0,\epsilon_1)$ together
with free summands. Margolis' proof implies also that every finitely-generated
$\mathbb{Z}$-graded $k[x,y]/x^2,xy,y^2$-module decomposes uniquely as a direct sum of
lightning flash modules. We note that $L(0,1,1)$ is the free $k[x,y]/x^2,xy,y^2$-module on one generator.

Hence we have:
\begin{eqnarray*}
\Rep(k[y]/y^2) & = & \mathbb{N}[\Sigma^{\pm 1}]\{ k,k[y]/y^2\} ,\\
\Rep(k[x,y]/x^2,xy,y^2) & = & \mathbb{N}[\Sigma^{\pm 1}]\{ \{ L(n,\epsilon_0,\epsilon_1): n\in\mathbb{N},\epsilon_0,\epsilon_1\in\{ 0,1\}\}\}, \\
\Rep(k[x]/x^2) & = & \mathbb{N}[\Sigma^{\pm 1}]\{ k,k[x]/x^2\}, \end{eqnarray*}
and the base-change maps are easily computed:
\[\xymatrix{ \Rep(k[y]/y^2) \ar[r] & \Rep(k[x,y]/x^2,xy,y^2)\ar[r] & \Rep(k[x]/x^2)  \\
 & L(n,0,0) \ar@{|->}[r] & \bigoplus_{i=0}^n \Sigma^{(|y|-|x|)i} k \\
 & L(n,1,0) \ar@{|->}[r] & \bigoplus_{i=0}^n \Sigma^{(|y|-|x|)i} k \\
 & L(n,0,1) \ar@{|->}[r] & \left(\bigoplus_{i=0}^n \Sigma^{(|y|-|x|)i} k\right) \oplus \Sigma^{(|y|-|x|)i} k[y]/y^2 \\
 & L(n,1,1) \ar@{|->}[r] & \left(\bigoplus_{i=0}^n \Sigma^{(|y|-|x|)i} k\right) \oplus \Sigma^{(|y|-|x|)i} k[y]/y^2 \\
k \ar@{|->}[r] & L(0,0,1) \ar@{|->}[r] & k[x]/x^2 \\
k[y]/y^2 \ar@{|->}[r] & L(0,1,1) \ar@{|->}[r] & k[x]/x^2. }\]
This verifies Corollary~\ref{graded applications corollary 2} directly, in this special case.

Every projective module over these algebras is free, hence on passing to the stable
representation monoids, we get the description:
\[\xymatrix{ \StabRep(k[y]/y^2)\ar[d]^{\cong} \ar[r] & \StabRep(k[x,y]/x^2,xy,y^2)\ar[d]^{\cong}\ar[r] & \StabRep(k[x]/x^2)\ar[d]^{\cong}  \\
\mathbb{N}[\Sigma^{\pm 1}]\{ k\} \ar[r] & \mathbb{N}[\Sigma^{\pm 1}]\{  \{ L(n,\epsilon_0,\epsilon_1): n\in\mathbb{N},\epsilon_0,\epsilon_1\in\{ 0,1\}\}\backslash \{ L(0,1,1)\}\} \ar[r] & \mathbb{N}[\Sigma^{\pm 1}]\{ k\} \\
 & L(n,0,0) \ar@{|->}[r] & \bigoplus_{i=0}^n \Sigma^{(|y|-|x|)i} k \\
 & L(n,1,0) \ar@{|->}[r] & \bigoplus_{i=0}^n \Sigma^{(|y|-|x|)i} k \\
 & L(n,0,1) \ar@{|->}[r] & \bigoplus_{i=0}^n \Sigma^{(|y|-|x|)i} k \\
 & L(n,1,1) \ar@{|->}[r] & \bigoplus_{i=0}^n \Sigma^{(|y|-|x|)i} k \\
k \ar@{|->}[r] & L(0,0,1) \ar@{|->}[r] & 0 . }\]
This verifies Corollary~\ref{graded applications corollary 3} directly, in this special case.
\end{example}

\bibliography{/home/asalch/texmf/tex/salch}{}
\bibliographystyle{plain}
\end{document}